\documentclass[12pt,reqno]{amsart}
\usepackage{CJKutf8}
\usepackage[margin=1in]{geometry}
\usepackage[T1]{fontenc}
\usepackage[utf8]{inputenc}
\usepackage{lmodern}
\usepackage{microtype}
\usepackage{amsmath,amssymb,mathtools,comment}
\usepackage{booktabs,array}
\usepackage{graphicx,tikz}
\usepackage{tabularx}
\usepackage{xcolor}
\usepackage{hyperref}
\usepackage{enumitem}
\usepackage{fancyhdr}
\usepackage{setspace}
\usepackage{titlesec}
\usepackage{url}
\graphicspath{{figures/}}
\newcommand{\red}{\color{red}}

\setlist{nosep,leftmargin=2em}
\titleformat{\section}{\Large\bfseries\sffamily}{\thesection.}{0.55em}{}
\titleformat{\subsection}{\large\bfseries\sffamily}{\thesubsection}{0.55em}{}
\titleformat{\subsubsection}{\normalsize\bfseries\sffamily}{\thesubsubsection}{0.55em}{}
\newcommand{\Ric}{\operatorname{Ric}}

\newcommand{\C}{\mathbb{C}}
\newcommand{\Q}{\mathbb{Q}}
\newcommand{\Z}{\mathbb{Z}}

\newcommand{\KeyIdea}[1]{\begin{quote}\small\itshape #1\end{quote}}

\newcommand{\optionalimage}[4]{%
  \IfFileExists{#1}{\includegraphics[width=#2,height=#3,keepaspectratio]{#1}}{%
    \fbox{\parbox[c][#3][c]{\dimexpr#2-2\fboxsep-2\fboxrule\relax}{%
      \centering\scriptsize #4\\JPG to come}}}}

\title[CY@50]{Calabi-Yau at 50:  mathematical discovery \\ in the age of AI}

\author{Yang-Hui He}
\address{London Institute for Mathematical Sciences, Royal Institution, London W1S 4BS, UK \&
Merton College, University of Oxford, Oxford OX1 4JD, United Kingdom}
\email{yh@lims.ac.uk}
\begin{document}

\begin{abstract}
Using the golden jubilee of Yau's proof of the Calabi conjecture as a spring-board, we discuss how the dataset of Calabi-Yau manifolds brought together mathematicians and physicists since the 1980s, and then further brought in the world of machine-learning and AI since 2017. 
We give a historical account on the growth of the field of ``AI-assisted mathematics'' in the last decade, and comment on what culminated in 2026 as perhaps the most dramatic year in the history of mathematics since Euclid.

Based on invited talk at {\it Calabi-Yau at 50: The Legacy and Future}, held at Tsinghua University in September, 2026.
\end{abstract}

\maketitle

\begin{quote}
\begin{flushright}
    \begin{CJK}{UTF8}{bsmi}
    春風滿座三千客，舊雨稱觴五十年\\
    (清) 王松《祝擎兄觀察五秩壽辰》
    \end{CJK}\footnote{
    By the caresses of spring's zephyr, three thousand do gather; to the tune of yesteryears' rain, we raise our glass to fifty autumns. \qquad -- Wang Song (1909)
    }
\end{flushright}
\end{quote}

\section{An anniversary and a Birth}
Fifty years mark a Golden Anniversary.
In the brevity of life and love, all too easily devoured by time, a half-century should indeed be cause to celebrate.
In 1977, in response to the problem formulated in the 1950s by Eugenio Calabi - who in 2023 passed away at the remarkable age of 100 - Shing-Tung Yau announced a solution, and gave its extended analytic proof in 1978, for which he was to be awarded the Fields Medal in 1982 \cite{Calabi1954,Calabi1957,Yau1977,Yau1978}.
The proof was not only a milestone in an ancient geometric tradition but also sealed a timely marriage of modern geometry and theoretical physics \cite{CHSW1985}.
Forty years later, the nuptials engendered, perhaps surprisingly, the first foray \cite{He2017} into AI-guided discovery in pure mathematics and theoretical physics, which grew in the last decade from a seedling to a forest which eclipses the entirety of science.

I am deeply honoured by the invitation of Professor Yau to give this talk on how Calabi-Yau manifolds came to be a nexus between mathematics, physics and AI, at {\it Calabi-Yau at 50: The Legacy and Future}, a conference held to celebrate the Golden Jubilee of the 1977 proof at Tsinghua/BIMSA.
I am confident that this happy union shall continue to bear fruit.

\subsection{Origins}
Our story perhaps begins with 1736, when Euler realized that a single integer, the genus $g \geq 0$, characterized the topology of (smooth, orientable, boundary-less) surfaces $\Sigma$; roughly it counts the number of ``holes'' as shown in Table \ref{tab:surface-trichotomy}.
Writing the eponymous combination
\begin{equation}\label{euler}
      \chi(\Sigma)=2-2g(\Sigma) \ ,
\end{equation}
the {\it Euler characteristic}, Gauss-Bonnet then links metric geometry with topology: total Gaussian curvature is $2\pi\chi(\Sigma)$. A quantity that  depends on the local bending of a surface is fixed by its global topology.

\begin{table}[tbp]
\centering

\label{tab:surface-trichotomy}
\renewcommand{\arraystretch}{1.18}
\begin{tabular}{@{}>{\centering\arraybackslash}p{0.31\textwidth}
                  >{\centering\arraybackslash}p{0.31\textwidth}
                  >{\centering\arraybackslash}p{0.31\textwidth}@{}}
\toprule
\textbf{Spherical} & \textbf{Flat} & \textbf{Hyperbolic}\\
\midrule
\optionalimage{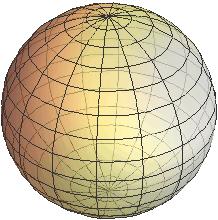}{0.17\textwidth}{0.3in}{sphere}
&
\optionalimage{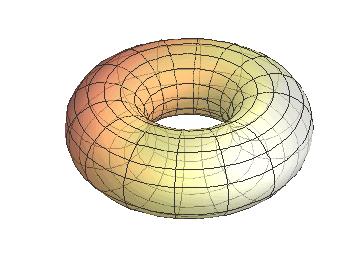}{0.17\textwidth}{0.3in}{torus}
&
\begin{tabular}{@{}c@{}}
\optionalimage{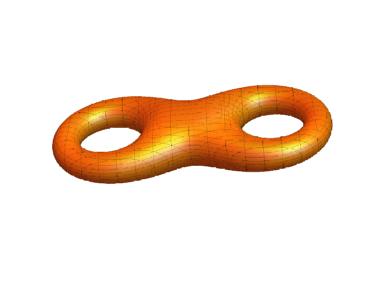}{0.075\textwidth}{0.7in}{$g=2$}\hspace{0.1em}
\optionalimage{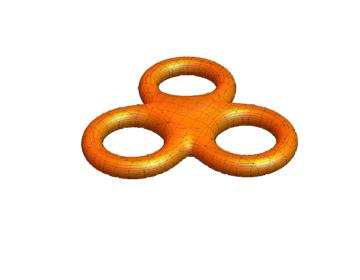}{0.075\textwidth}{0.48in}{$g=3$}\hspace{0.1em}
\optionalimage{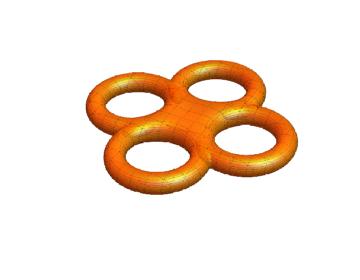}{0.075\textwidth}{0.48in}{$g=4$}
\end{tabular}
\\
$g(\Sigma)=0$ & $g(\Sigma)=1$ & $g(\Sigma)>1$\\
$\chi(\Sigma)=2$ & $\chi(\Sigma)=0$ & $\chi(\Sigma)<0$\\
positive curvature & zero curvature & negative curvature\\
sphere & Ricci-flat torus & higher-genus surface\\
\bottomrule
\end{tabular}
\caption{The genus trichotomy for compact orientable surfaces.}
\end{table}

Eq.~\ref{euler} in conjunction with Gauss (q.v.~Table \ref{tab:surface-trichotomy}) thus gives a sign trichotomy: the sphere has positive curvature, the torus has zero curvature, and more handles imply negative curvature.

Riemann then adds a complex structure: real two-dimensional surface becomes a Riemann surface, a one-dimensional complex manifold. The above trichotomy reappears in a new vocabulary.
The sphere is the projective line $\mathbb{C}\mathbb{P}^1$ and is the simplest Fano variety.
The torus is Ricci flat and is the simplest Calabi-Yau. 
The remaining, hyperbolic, negatively curved cases, are so-called general type.

In the 20th century, more branches of mathematics enter the scene.
The integral of Gaussian curvature is seen, in modern algebraic geometry, as an inner product between the first Chern class $c_1$ of the tangent bundle $T\Sigma$ and the class $[\Sigma]$ of the curve itself.
The Euler characteristic is seen as an alternating sum (index) of (co-)homological objects. 
Number-theoretically,  Deligne, Dwork, Grothendieck, and Weil
have shown that the zeta-function $Z_\Sigma(s)$- the counting function of $\Sigma$ over finite fields - also sees the topology in its rational form. We summarize these ideas in Table \ref{tab:surface-lenses}, marking in red the recurrent motif of the Euler characteristic.

\begin{table}
\centering
\label{tab:surface-lenses}
\small
\begin{tabular}{c|c}
\toprule
Topology &
$\displaystyle {\red \chi(\Sigma)}=2-2g(\Sigma)$\\
Algebraic geometry &
$\displaystyle \langle c_1(T\Sigma),[\Sigma]\rangle={\red \chi(\Sigma)}$\\
Differential geometry &
Gauss-Bonnet: $\displaystyle \frac{1}{2\pi}\int_\Sigma R\,dA={\red \chi(\Sigma)}$.\\
Homological algebra &
Index: $\displaystyle {\red \chi(\Sigma)}=\sum_{i=0}^{2}(-1)^i \dim H^i(\Sigma;\mathbb{Q})$.\\
Number theory &
$\displaystyle Z_\Sigma(s)=\exp\!\left(\sum_{m\geq1}\frac{\#\Sigma(\mathbb{F}_{q^m})}{m}q^{-ms}\right) \leadsto
Z_\Sigma(1-s)=q^{(1/2-s){\red \chi(\Sigma)}}Z_\Sigma(s)$\\
\bottomrule
\end{tabular}
\smallskip

\caption{The Euler characteristic $\chi$ seen through five mathematical lenses.}
\end{table}

Generalizing Table \ref{tab:surface-lenses} to beyond complex dimension 1 is a unifying theme of contemporary mathematics, inspiring, inter alia, the famous minimal model programme and the Langlands programme.
While this is a momentous task in general, at least for a subclass of complex manifolds, viz., K\"ahler, there has been steady progress.
This is why our story is so relevant: it is one which threads a narrative over three centuries.

\subsection{Calabi and Yau}
A K\"ahler manifold is a complex manifold which has a metric compatible with its complex structure in a way that allows geometry to be encoded by a closed real $(1,1)$-form, the Kahler form $\omega$. Locally, in holomorphic coordinates $z^1,\ldots,z^n$, the metric coefficients come from a potential $K$:
\[
  g_{i\bar j}=\frac{\partial^2 K}{\partial z^i\partial\bar z^j},
  \qquad
  \omega=\sqrt{-1}\,g_{i\bar j}\,dz^i\wedge d\bar z^j.
\]
The potential is not unique: adding the real part of a holomorphic function changes no metric. Yet this local freedom conceals a global constraint. The cohomology class $[\omega]$ is fixed as one moves among metrics in the same Kahler class. The question becomes: within the allowed class, can one choose a metric whose Ricci curvature realizes a prescribed representative of the first Chern class?

This is Calabi's question \cite{Calabi1954,Calabi1957}. 
In dimension one, the sign of the canonical curvature organizes the answer beautifully; in higher dimensions, Chern classes take over the organizing role. The surface diagram is not a toy to be discarded when we reach threefolds. It is a small model of the sort of correspondence we would like to generalize.

For a physicist, the phrase ``Ricci-flat'' may ring a different bell. General relativity teaches that curvature is linked to energy and matter; a vacuum Einstein metric has vanishing Ricci tensor. On a Kahler manifold, Ricci-flatness is more constrained than it sounds because the complex structure must also be respected. A Calabi-Yau space therefore offers a meeting point for topology, complex analysis, nonlinear partial differential equations, and mathematical physics.

The Calabi conjecture, in its standard formulation, concerns a compact Kahler manifold $(M,\omega)$ and a smooth real $(1,1)$-form $\rho$ representing the real first Chern class. It asks whether there is a unique Kahler metric $\widetilde\omega$ in the same cohomology class whose Ricci form is $\rho$. Calabi established uniqueness and made progress toward existence; Yau proved the existence statement, resolving the conjecture \cite{Calabi1954,Calabi1957,Yau1977,Yau1978}.

\subsection{The hyphen}
Now the hyphen. In the 1980s, in a seemingly unrelated development in theoretical physics over the decades since quantum mechanics, string theory gave complex algebraic geometry an unexpected address in physics. A ten-dimensional theory needed a way to leave four large spacetime dimensions and curl up six internal ones. In the supersymmetric compactifications of the heterotic superstring by Candelas-Horowitz-Strominger-Witten, the internal space had to satisfy stringent holonomy conditions \cite{CHSW1985}.
It so happened that Strominger and Yau were next door at the IAS and their conversations brought Ricci-flatness and K\"ahlerity, and indeed complex algebraic geometry, into physics.
It is curious that it was the physicists who coined the term {\it Calabi-Yau}, to describe (compact smooth) Ricci-flat K\"ahler manifolds.

\subsection{The bestiary, the database, and the Neural Network}
The decades between the 1980s and the turn of this century saw intense activity resulting from the birth of the union between physics and mathematics: the explicit construction of large data-bases of Calabi-Yau manifolds.
To the distinguished company here I am, of course, preaching to the choir, so I shall skip the details on the fruits of these many years, and leave the students to the classic ``bestiary'' by H\"ubsch \cite{Hubsch1992}, as well as a contemporary homage thereto \cite{He2021Landscape}.

I realize that I am supposed to talk more about the AI aspects of theoretical discovery, so I will now move to the main subject of my lecture and hope the audience would see the above as a long introduction to pay my respects to a happy 50 years of Calabi-Yau.
Though I referred to a birth in the above, perhaps a birth more pertinent here is 2017, the fortieth anniversary.

As mentioned, the Bestiary \cite{Hubsch1992} did not only document the Calabi-Yau union between mathematics and physics, but it also accompanied the establishment of large databases of Calabi-Yau data, numbering in the billions.
I spent much of the earlier part of my career studying such data, both the rich geometry in her own right, and to mine it to try to find universes akin to our own via string compactifications. 

In the 2010s, students and postdocs in mathematics and physics, unlike the generation who moved into finance before the 2008 global crisis, would generically move onto {\it data science} and related {\it machine learning}.
I had no idea what these words meant but was curious to find out.
In 2017 - which happened to be my son's birth - sleepless nights from a young daughter and an infant son induced me to take online tutorials on machine learning.
Immediately after mastering everyone's first baby-step of writing a relatively shallow neural network in Python to classify hand-written digits, I tried the perhaps bizarre thing of passing Calabi-Yau data into the input layer.
This, in retrospect and to my knowledge, became the first intersection between the world of pure mathematics and the world of AI.

In brief, a Calabi--Yau threefold comes with a Hodge diamond, the (anti-)holomorphic refinement of the Betti numbers. The two non-trivial numbers therein, $h^{1,1}$ and $h^{2,1}$, are especially familiar: in the standard setting they count Kahler and complex-structure moduli, and mirror symmetry exchanges them.
Moreover, they determine the Euler characteristic as $\chi = 2(h^{1,1} - h^{2,1})$.
A typical exercise in algebraic geometry and related string compactification, is to compute, via long exact sequence chasing \footnote{The Euler characteristic itself is somewhat more easily determined by the index theorems that generalizes Gauss-Bonnet by an integral. The individual (co-)homological terms are much more difficult to obtain.} - the sort of rite of passage suffering for the beginning PhD student - say $h^{2,1}$ from a given configuration for the manifold.
For instance,
\begin{equation}    
X={ \left(\arraycolsep=1.4pt\def\arraystretch{0}
\begin{array}{cccccccc}
 1 & 1 & 0 & 0 & 0 & 0 & 0 & 0 \\
 1 & 0 & 1 & 0 & 0 & 0 & 0 & 0 \\
 0 & 0 & 0 & 1 & 0 & 1 & 0 & 0 \\
 0 & 0 & 0 & 0 & 1 & 0 & 1 & 0 \\
 0 & 0 & 0 & 0 & 0 & 0 & 2 & 0 \\
 0 & 1 & 1 & 0 & 0 & 0 & 0 & 1 \\
 1 & 0 & 0 & 0 & 0 & 1 & 1 & 0 \\
 0 & 0 & 0 & 1 & 1 & 0 & 0 & 1 \\
\end{array}
\right)}
\longrightarrow
\mbox{$H^*$ induced from Euler sequence }
\longrightarrow
h^{2,1}(X) = 22
\label{CICY2000}
\end{equation}

In the above example, I have taken a so-called CICY $X$ (complete intersection Calabi-Yau threefold), where the columns record the multi-degrees of $X$ embedded into products of projective spaces.

\subsection{Mathematics as Image-Processing}
The details are important, but every Calabi-Yau manifold can be interpreted as a configuration matrix like Eq.~\ref{CICY2000}. Indeed, any algebraic variety, whether defined combinatorially or algebraically, can be thus represented.

My excursions into machine learning immediately made me rethink such configurations as pixelated images \cite{He2017}, then Eq.~\ref{CICY2000} essentially becomes an image labelling problem:
\begin{equation}\label{NNhodge}
\begin{array}{c}
\includegraphics[width=0.2\linewidth]{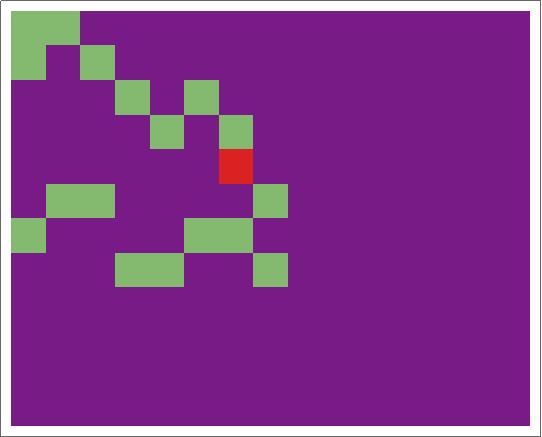}
\end{array}
\longrightarrow
\begin{array}{c}
\begin{tikzpicture}[
    neuron/.style={circle, draw, minimum size=1cm, inner sep=0pt, thick},
    input/.style={neuron, fill=blue!20},
    hidden/.style={neuron, fill=green!20},
    output/.style={neuron, fill=red!20}
]

\def\inputnum{3}
\def\hiddennum{4}
\def\outputnum{1} 

\foreach \i in {1,...,\inputnum} {
    \node[input] (I-\i) at (0, -\i) {$x_\i$};
}

\foreach \j in {1,...,\hiddennum} {
    \node[hidden] (H-\j) at (3, -\j+0.5) {$h_\j$};
}

\foreach \k in {1,...,\outputnum} {
    \node[output] (O-\k) at (6, -2) {$y_\k$};
}

\foreach \i in {1,...,\inputnum} {
    \foreach \j in {1,...,\hiddennum} {
        \draw[->, draw=gray] (I-\i) -- (H-\j);
    }
}

\foreach \j in {1,...,\hiddennum} {
    \foreach \k in {1,...,\outputnum} {
        \draw[->, draw=gray] (H-\j) -- (O-\k);
    }
}

\node at (0, 0.5) {\textbf{Input Layer}};
\node at (3, 0.5) {\textbf{Hidden Layer}};
\node at (6, 0.5) {\textbf{Output Layer}};
\end{tikzpicture}
\end{array}
\longrightarrow
22
\end{equation}
For the image to the left which represents $X$, I have used a simple representation that 0 is purple, 1, green and 2 red.
Any other choice of data representation of $X$ is, of course, equally valid and the performance of the machine learning certainly would depend on the choice; I was lucky that my naive choice was good enough.

Furthermore, I was fortunate that I had millions of such labelled data, benefitting from the work of the community from the 1980s to 2000s, and in addition that GPUs introduced to the personal laptop in the 2010s were powerful enough to handle such classification problems in a matter of minutes to hours.
Accuracies, after cross-validation, quickly reached the nineties.
I was rather surprised: how could a neural network which knows nothing about algebraic geometry, get so good at guessing maths?

Fanaticism ensued from enthusiasm, I started passing every type data from every branch of mathematics I could understand into neural networks and other machine learning algorithms. Since 2017, I made scores of new academic friends, writing some 100 papers on AI-guided mathematics and theoretical physics in due course. After all, the only thing more fun than doing research is making friends whilst doing it.
There is clearly not enough space here for me to thank each of my wonderful students and collaborators by name; I leave that to a collected Google-doc \cite{thanks} that acknowledges the people and the papers.

At some level, one could even argue that all of mathematics within a given order of logic, whether it be sentences of proof and symbols, or data that represents everything from sequences of prime numbers to polynomial coefficients of algebraic varieties, can be tensorized and fed into a Turing machine.
Now, instead of an idealized Turing machine, we have neural networks that approximate.
I joked back then that algebraic geometry, if not all of mathematics, is an image processing problem, and I was pleasantly surprised that such a statement was not met with hostility.

Now, 2017 was an auspicious year, in addition to the two stories mentioned here, the fortieth anniversary of Yau's proof and the birth of my son that led my AI experimentation with mathematics, there was the android Sophia becoming the first non-human to obtain a passport \cite{bbc2017sophia}, there was the publication of the landmark paper \cite{vaswani2017attention} on Transformer architecture which underpins all of the modern LLMs from ChatGPT to Claude.
Indeed, my colleagues were independently exploring the string theory landscape with machine learning and AI algorithms \cite{KreflSeong2017,Ruehle2017,Carifio2017}.
An immediate by-product of this coincidence was the conference \verb|String_Data|, which has since become the first and still-running annual conference on what is now called ``AI for theoretical physics'' \cite{stringdata}.
Not long after, the annual series DANGER, starting in 2021, which stands for Data in analysis, geometry and representation theory, became one dedicated to ``AI for pure mathematics'' \cite{danger}.

I hope the audience can forgive my shameless self-promotion in referring them to the textbook \cite{He2021Landscape}, the editorial \cite{He2023ML} and the tutorial \cite{constantin2026machine}, for the development over the past decade.
In the ensuing, I shall largely adhere to a trinity of how AI can assist mathematical discovery in a review I wrote in 2024 \cite{He2024Triumvirate}, as well as a more recent commentary with my colleagues at the London Institute earlier this year \cite{burtsev2026how}. These, in the rapidity of what has transpired over the last few months, are already rather outdated. I shall finish with commentary on the so-called ``MathApocalypse'' of October 6th, 2026.

\section{A Triumvirate of Mathematical Discovery}
Where, then, are we? It has been almost ten years since the sleepless nights of 2017, and the habit of applying machine learning to pure mathematics and theoretical physics has grown from an eccentricity into an industry.

The phrase ``AI'' means many things to many people, and I shall not be so foolhardy as to define it. To some it simply means language models; to others, neural networks or the Transformer architecture; to others still, anything a computer does that impresses them. What I can say is that the words first came together in the proposal for the summer meeting at Dartmouth in 1956 \cite{Dartmouth1955}, a date we shall meet again presently. There is an old adage among computer scientists, usually attributed to Larry Tesler, that AI is whatever has not been done yet. 
For decades this was a comfortable definition, since it guaranteed that AI never arrived. It has now very much arrived.

Taking mathematics in the generous sense which includes theoretical physics, I have found it helpful to organise the ways in which AI can assist discovery into three directions \cite{He2021MLMS,He2024Triumvirate}, summarised in Table~\ref{tab:triumvirate}. The words are used loosely, and the boundaries between them, as we shall see, have all but dissolved in the last year. Nevertheless, they remain a serviceable map.

\begin{table}[t]
\centering
\renewcommand{\arraystretch}{1.25}
\small
\begin{tabular}{@{}p{0.2\textwidth}p{0.28\textwidth}p{0.42\textwidth}@{}}
\toprule
\textbf{Direction} & \textbf{Mathematics as \ldots} & \textbf{What the machine does}\\
\midrule
Bottom-up & a formal logical system & building proposition symbol by symbol, from axioms\\
Top-down & a creative, intuitive art & experimental mathematics: finds patterns, conjecture and new directions\\
Meta-mathematics & a language & LLMs on maths with or without ``understanding''\\
\bottomrule
\end{tabular}
\smallskip
\caption{The triumvirate of AI-assisted mathematical discovery, following \cite{He2024Triumvirate}.}
\label{tab:triumvirate}
\end{table}

When I wrote the review \cite{He2024Triumvirate} in 2024, I thought of the three as roughly comparable in promise, with the third the least mature. That review is now grossly outdated, and so, I fear, is the commentary \cite{burtsev2026how} which my colleagues and I wrote only months ago. It is an extremely exciting time, and we are very lucky to be alive in it. I shall argue, with evidence, that 2026 is a landmark in the history of mathematics, perhaps the most significant since Euclid. But let me first take the three directions in turn.

\subsection{Bottom-up: Axioms and Formalization}\label{sec:bottomup}
The first direction is what Terence Tao has come to call the {\it proof co-pilot} \cite{Tao2025MAP}.
Bottom-up mathematics is mathematics built from its foundations: one begins with axioms and definitions and seeks a path, each step licensed by the rules of inference, towards true statements. This is a programme which goes back to Euclid, and which passed through Leibniz and Frege on its way to the twentieth century's answer to the {\it Elements}: the {\it Principia Mathematica} of Whitehead and Russell \cite{WhiteheadRussell}. It is a book which, unless one is a logician, nobody has read, in part because its ambition to axiomatise all of mathematics was killed within twenty years by G\"odel's incompleteness \cite{Godel1931} and by the undecidability results of Church and Turing \cite{Church1936,Turing1936}.

And yet. I quote my very learned friend Prof~Minhyong Kim: ``The practicing mathematician hardly ever worries about G\"odel.'' If one had to ask of every lemma whether it were decidable, one would become very depressed very quickly. Yes, there are undecidable statements; there are also infinite numbers of interesting statements which are perfectly provable, and it is on these that we spend most of our days.

So, G\"odel notwithstanding, what used to be called automated theorem proving (ATP) has a long and honourable tradition. It goes back, to 1956, the very year of the Dartmouth meeting \cite{Dartmouth1955} which induced the words ``Artificial Intelligence'' to be first written in print, when Newell, Simon and Shaw devised their Logic Theory Machine, which proceeded to prove a number of the early theorems of the {\it Principia} \cite{NewellSimon1956}. The 1970s brought the type theories of Martin-L\"of and, later, Coquand, upon which the modern proof assistants rest; the 2000s brought Voevodsky's univalent foundations and homotopy type theory \cite{HoTT2013}. By then the Coq system had dispatched the four-colour theorem \cite{Gonthier2008} and, in 2012, the Feit-Thompson odd-order theorem \cite{Gonthier2013}, whilst Hales and his collaborators formally verified the Kepler conjecture \cite{Hales2017}. These were heroic efforts. They were also, it is fair to say, regarded by most working mathematicians as somebody else's business.

\subsubsection{Mathlib}
The twenty-first century's answer to the {\it Principia}, in the lineage of Hilbert's formalism and its defiant ``we must know, we will know'', is Mathlib \cite{mathlib2020}: the community library of the Lean proof assistant, which was launched in 2013 by Leonardo de Moura \cite{deMoura2015}. 
The principle is disarmingly simple. One rewrites one's mathematics, every definition and every line of every proof, in Lean code, and if the code faithfully represents the problem and it compiles, then the theorem is correct. There are other such languages, Coq and Isabelle among them, but Lean is the one which mathematicians have come to love. Every one of my devices now carries a Lean sticker.

That love owes much to one of the torchbearers of the community, my colleague and friend Kevin Buzzard, who is unmistakable at any conference for his very brightly coloured trousers. Through his Xena Project, Buzzard recruited undergraduates at Imperial College, a couple of hundred every year, to formalise literally all of undergraduate mathematics \cite{Buzzard2022ICM}. This was human labour: thousands of people over the years, typing in every statement and every proof of the standard curriculum. In 2023 they pressed the return key, as it were, and it compiled. Thank goodness. One shudders to imagine what would have followed had the theorems we teach our first-years contained mistakes.

Proof assistants, as I said, are as old as AI itself. But something changed in 2019, when Buzzard, Commelin and Massot formalised perfectoid spaces \cite{BCM2020}. Kevin cheerfully calls this a marketing strategy, and it was a brilliant one: perfectoid spaces are the objects for which Peter Scholze was awarded the Fields Medal \cite{Scholze2012}. For the first time, contemporary mathematics of the highest order, and not merely textbook material, had been made to compile. Overnight the community shifted from computer scientists and logicians to mainstream mathematicians, and Mathlib became the thing to do. 
The point I wish to emphasise is one of trust. It is an uncomfortable fact that there are, by now, papers in the very best journals which are simply wrong. We have come to trust a proof that compiles against Mathlib at least as much as we trust the most severe peer review. It is therefore of the utmost importance that, sooner rather than later, everything should go in.

\subsubsection{The Lean community}
The growth since has been remarkable, and a few dates may serve as signposts.
In 2023 the Lean Focused Research Organization was founded. One sees immediately where this must go: if formalisation is to keep pace with mathematics, one cannot rely on armies of undergraduates; one must train a language model to translate \LaTeX\ into Lean. This is {\it auto-formalization} \cite{Szegedy2020,Lu2024}, and it began to work in earnest this year, with the appearance of, for instance, \cite{CFSG}. It is one of the main reasons why everything else in this article has happened so fast.

In 2024, Buzzard disappeared into his attic, supported by a multi-year EPSRC grant, to formalise the proof of Fermat's Last Theorem \cite{FLTproject}, taking for granted, for the time being, the mathematics known up to about the 1980s. Here one meets a charming feature of Lean. If one does not yet have a proof, one simply writes the word \verb|sorry|, and carries on. The Riemann Hypothesis, for example, has a statement in Lean whose proof currently reads \verb|sorry|; yet one may assume it and proceed to formalise, quite rigorously, every result which is conditional upon it.

Then, in the week before this lecture, came the news that Claude had auto-completed the entire Fermat project in eleven days, producing some fourteen million lines of verified Lean code.
Kevin, I am glad to report, is still not depressed. The machine, he points out, took a somewhat different route through the literature from the one he and his team are following, so that there is work yet to be done. Still, one sees how fast this is going: a project designed to occupy a team of experts for the better part of a decade was, in its essentials, overtaken in days. Ironically, only at the beginning of 2025, Kevin still insisted that he had never used any AI in Mathlib.

For the theoretical physicists, I warmly recommend Tooby-Smith's PhysLean \cite{PhysLean}, which is Lean for physics. The \verb|sorry| mechanism is a godsend here. Define a quantum field theory: \verb|sorry|, we do not have that definition yet. But assume that it exists, and compute a Feynman diagram. Many of the projects I now suggest to younger colleagues are of this kind: to formalise the papers of theoretical physics, which is among the most depleted regions of the library. Mathlib has modular forms. It has perfectoid spaces. Yet as of January of this year it did not have a definition of the Bessel function, which is our bread and butter, simply because the community had not been there to type it in.

I have the honour of helping to run the (other) LLM, which stands for the ``London Lean Monthly'', where Buzzard, Mehta, Hill and I gather the Imperial, UCL and London Institute communities once a month to take stock of formalisation.
Finally, and this will matter in the sequel, DeepMind's AlphaProof \cite{AlphaProof2025}, unlike its geometric cousins AlphaGeometry and AlphaGeometry2 \cite{AlphaGeometry2024}, was built around Lean from the outset. Seen in this light, bottom-up mathematics is a reinforcement-learning environment: one has a statement, and one searches for a pathway of correct intermediate statements until the theorem is reached, with the compiler as an incorruptible referee.

\subsection{Top-down: Intuition and Data}
So much for building from the ground up. But mathematics, as every practitioner knows, is also a top-down process; it is experimental. With a small nod to Brouwer and the intuitionists, for whom mathematics is a creation of the mind, one must concede that we rarely discover theorems by deduction from axioms. We guess, on the strength of intuition and experience, and justify afterwards. 
In parallel, theoretical physics has {\it never} been a bottom-up process. I love the remark of the great Arnol'd that ``mathematics is the part of physics where experiments are cheap'' \cite{Arnold1998}.

Cheap, and also perfect. It was most fitting that Mike Douglas spoke at this conference shortly before I did, for it is he who has popularised a term which I adore: {\it Platonic data}. Mathematical data has no experimental error and no noise. It is simply out there, exact, awaiting inspection. What I was doing back in 2017, was the beginning of AI exploration of Platonic Data.

\subsubsection{Intuition, experience, experimentation}\label{sec:topdown}
The practice is old. In the 1790s, the teenage Gauss stared at tables of primes and at the staircase $\pi(x):=\#\{p\leq x\}$, without a computer and well before complex analysis, and immediately saw that it was approximately $x/\log x$.
He recorded the observation in a letter to Encke in 1849, and one then had to wait another half-century for Hadamard and de la Vall\'ee Poussin, in 1896, to prove what is now the Prime Number Theorem,
\begin{equation}
\pi(x)\sim \frac{x}{\log x} \ .
\end{equation}
They had to wait because the method of proof required complex analysis, and complex analysis had to wait for Cauchy and Riemann. The brain of Gauss was, I like to say, the best neural network of the nineteenth century. I exaggerate only a little: it is no accident that the name of Gauss is attached to least squares and to the bell curve, for the fitting of curves to data was very much his craft. 

In the twentieth century, the Birch-Swinnerton-Dyer (BSD) conjecture was already computer-aided conjecture formulation on Platonic data. In the early 1960s, Birch and Swinnerton-Dyer went down to the computer in Cambridge and plotted, for elliptic curves, the rank $r$ against the point counts $N_p$ modulo primes \cite{BSD1963}. Out came one of the greatest problems in mathematics, which remains open; for how much longer, I no longer dare to say.

It bears emphasis that a third of the remaining Millennium Prize Problems came from experimentation: BSD from those plots, and the Riemann Hypothesis from Riemann computing by hand the first few zeros and deciding that they looked fine. I used to write ``two of the remaining six'' on this slide. For reasons which will become apparent in Section~\ref{s:2026}, I am no longer sure whether the denominator should be six or five.

Today Platonic Data is abundant. There is of the order of a terabyte of it freely available online, for the wonderful thing about the mathematics community is that we give everything away: GAP for finite groups \cite{GAP}, the Graded Ring Database and the Calabi-Yau databases for algebraic geometry \cite{GRDB,He2021Landscape}, the LMFDB for number theory \cite{LMFDB}, SnapPy for knots and three-manifolds \cite{SnapPy}, and so on. I can only encourage the student to pick a favourite area and start mining.

\subsubsection{The Birch Test}
In the autumn of 2023 a number of us, among them Buzzard, Cardoso, Klemm and Nampuri, organised a programme at the Isaac Newton Institute in Cambridge on black holes, number theory and holographic quantum information \cite{INI2023}. The title was meant to end with ``and machine learning'', but this exceeded website's permitted number of characters and was duly chopped off. The purpose was to bring friends from the maths, physics, data-science, and AI communities into one room and let them talk. There we watched a recorded lecture by Bryan Birch who, I am delighted to say, is well into his nineties and going strong; mathematicians do tend to be long-lived. 
Calabi lived to 100 and I hope Prof. Yau will have many more happy returns.

At the conference, since ChatGPT had, by common consent, passed the Turing test in 2022, we thought it expedient to propose, half-jokingly, a new benchmark for AI-driven mathematical discovery. We called it the {\it Birch Test}, and Burtsev and I reported it in {\it Nature} at the start of 2024 \cite{HeBurtsev2024}. It asks that a discovery be ``A, I, N @ AI'' :
\begin{itemize}
\item[(A)] {\bf Automatic}: it must be generated by the AI, without humans mucking about in the process;
\item[(I)] {\bf Interpretable}: it must be concrete enough to be stated as a conjecture;
\item[(N)] {\bf Non-trivial}: it must be hard and interesting enough for the community to work on it.
\end{itemize}
Most importantly, it must make Birch happy.

The middle criterion deserves comment. Something like \eqref{NNhodge} fails it immediately. The network predicted Hodge numbers to splendid accuracy, but what was it telling us by way of new theorems or new ideas? I did not know. Accuracy without interpretation is an oracle, and an oracle is not mathematics.

As of now, there hasn't been a major conjecture raised by AI that passes the Birch Test.
There have been some near misses. The beautiful work of Davies et al.\ at DeepMind found new relations among knot invariants by examining saliency in networks trained on SnapPy data \cite{Davies2021}; it fails (N) in the sense that the authors were able to prove the result themselves. The discovery of unstable singularities in fluid equations \cite{Wang2025} fails (A), since there was a good deal of expert human steering.

There is one close to my heart. With Kyu-Hwan Lee, Tom Oliver and Alexey Pozdnyakov, we found a new and wholly unexpected pattern in the $L$-functions of elliptic curves, which generalises Chebyshev's bias in the distribution of primes \cite{He2025Murmuration}. Averaged over curves of given rank, the Frobenius traces oscillate in a manner which, when plotted, resembles  a flock of starlings at dusk. This has come to be known as the {\it murmuration phenomenon}. It fails (A), because we were mucking about with the discovery process. But it certainly passes (N): it is, to my knowledge, the first conjecture raised with the help of AI which remains open and upon which serious number theorists continue to work \cite{Chiou2024}. The mathematics has long since risen above my expertise in number theory. I think Birch was almost happy. My proudest moment, however, was when my friend, the famous composer Cheryl Frances-Hoad, wrote a concerto for saxophone and strings whose second movement is entitled {\it Elliptic murmurations} \cite{FrancesHoad2024}; it was performed at the World Science Festival in Australia this March for which I was honoured to be a part \cite{ANYO}.

\subsection{Meta-mathematics: LLMs on Maths}
This brings me to the third direction and to the heart of the matter.
When I was asked to write a review  \cite{He2024Triumvirate}, two years ago, I think it is fair to say that most people treated large language models (LLMs) with more contempt than admiration. They were fine, they were amusing, but they made every kind of mistake and they hallucinated all the time. I allowed that perhaps, one day in the foreseeable future, they might fundamentally change the landscape of mathematical research. What is truly remarkable is that this third direction has not only caught up with the other two but has swallowed them. Today, LLMs have become synonymous with AI.

The pace may be judged from Table~\ref{tab:rapid}. A word on its first line. The imitation game was proposed by Turing in 1950 \cite{Turing1950},  when it was inconceivable that anything should pass it. It is, however, not a very rigorous test, and we passed it in 2022 not because we had designed an intelligent machine but because of statistics: we had extracted enough data for a machine to fool a human in conversation. Then, in 2025 and 2026, everything mushroomed, and private companies entered in force. It is worth understanding why. It is not, I am sorry to say, because they find mathematics interesting. It is because they wish to sell a product whose selling point is {\it verifiability}, and mathematics is the cleanest benchmark there is for whether a machine is correct.

\begin{table}[t]
\centering
\renewcommand{\arraystretch}{1.2}
\small
\begin{tabular}{@{}lp{0.8\textwidth}@{}}
\toprule
2022 & ChatGPT passes the Turing test\\
2023 & DeepMind's FunSearch \cite{FunSearch2024}; Meta AI's Llama: language models begin to do mathematics\\
2024 & DeepMind's AlphaGeometry and AlphaGeometry2 \cite{AlphaGeometry2024} solve Olympiad geometry\\
2024--5 & AlphaProof \cite{AlphaProof2025} incorporates Lean\\
2025 & AlphaEvolve \cite{AlphaEvolve2025}, tailored to optimisation problems\\
2025--6 & Mushrooming of companies: Axiom; Aristotle at Harmonic; Gauss at Math Inc.; Aletheia at DeepMind; Archon; Danus at Rethlas (all incorporating Lean); ProofCouncil (incorporating an agentic harness)\\
\bottomrule
\end{tabular}
\smallskip
\caption{Rapid progress of LLMs for Maths}
\label{tab:rapid}
\end{table}

\section{Vibe Maths}
Why should a language model be any good at mathematics at all?
At school we are taught the standard methods of proof: by contradiction, by induction, by syllogism, etc.
But there is a wonderful essay by Thurston, who shared the Fields podium with Yau in 1982, entitled ``On proof and progress in mathematics'' \cite{Thurston1994}. He points out that although mathematics presents itself as a bottom-up edifice, much of what mathematicians actually do is top-down. One has a gut feeling as to whether something is true. Physicists live by this; we have a feeling for what is correct long before, and sometimes entirely without, a proof.

In practice, as every seminar-goer knows, the methods of proof are rather more various. This is a joke, obviously, but we prove by intimidation, by authority, by citation, by obfuscation and, most importantly, we prove by {\it vibe}. That is the gut feeling, and it is the word which has caught on. LLMs certainly learnt proof by citation and proof by authority, and learnt them all too well.
The question for meta-mathematics is therefore this: how do we make proof by vibe reliable?

\subsection{Outdated jokes}
In July of 2025, I was telling audiences that AI was approximately as good as the twelve-year-old Terry Tao. The reference was to the International Mathematical Olympiad, at which Tao famously collected bronze, silver and gold before he was properly a teenager. DeepMind's systems scored 28 points out of 42 at the IMO in 2024, the standard of a silver medal, and 35 points in 2025, which is gold.
I used to add that when it reaches the sixteen-year-old Terry, we should start to worry.

There is an earlier joke, by Christian Szegedy \cite{Szegedy2020}. In 2020 he made an extraordinary extrapolation: computers surpassed humans at chess in the 1990s and at Go in the 2010s; 
they should therefore surpass humans at proving theorems by 2030. Our reaction at the time was that this was meant as a joke; after all there were only two data points.

What is remarkable is that Szegedy's estimate now appears to have been conservative. Both jokes are outdated.

\subsection{A phase transition}
The general feeling of the community is that the last eight or nine months have seen a phase transition in the so-called reasoning ability of AI. Whatever one means by {\it reasoning} I leave to the philosophers; let us consider the practical side. What is different between the systems of today and the ChatGPT of three years ago? I can identify a few ingredients:

\begin{description}[leftmargin=2em,style=nextline,itemsep=0.4em]
\item[Agentic AI]
One no longer converses with a single model. Many agents, built upon LLMs, work together and check one another: there is the prover, the proof-checker, and the checker of the checker.

\item[Harnessing]
Also called scaffolding: the ability of the model to call external programmes, such as Sage, Python and, most importantly, Lean. Three years ago, if one asked ChatGPT for the seventh digit of $23/145$, it would free-associate, with a one-in-ten chance of success. Today it simply hands that part of the computation to something which can do arithmetic.

\item[RLVR over RLHF]
Reinforcement learning with verifiable rewards, in place of human feedback. Instead of human beings labelling what is good and what is bad, the reward comes from something which can be checked. In mathematics the supreme verifier is the Lean compiler. This is precisely where the bottom-up direction of the triumvirate is absorbed into the third.

\item[Test-time compute]
The ``P'' in GPT stands for ``pre-trained'': one trains the network, freezes it, and what one then interrogates is a rigid object. Now, ever more of the computational resources are spent at the time of the question rather than at the time of training. The machine, in the vernacular, thinks for longer.

\item[GRPO]
Group Relative Policy Optimisation \cite{GRPO2024} is one of the mechanisms by which the system comes to grade its own intermediate steps, and not only its final answer. If one wants a working definition of reasoning, that is not a bad one.
\end{description}

The combination of these, all coming to maturity at once, has turned what used to be a joke about hallucinating chatbots into something highly non-trivial. Proof by vibe, constrained by an agentic harness and refereed by a compiler, turns out to be a rather good description of how mathematics is done.  {\it Vibe maths} by AI is getting as good as vibe maths, in the sense of Thurston, by humans. And in addition, AI can check using Lean, which is better than humans can.

\subsection{2026: A Turning Point in History}\label{s:2026}
What, concretely, has been achieved? Table~\ref{tab:score} is the scoreboard as it stood on the day of the lecture in September; the numbers will certainly have changed by the time this is read.
The one I know best is FrontierMath, and since it has served as the best thermometer of the phase transition, let me tell its story.

\begin{table}[t]
\centering
\renewcommand{\arraystretch}{1.2}
\small
\begin{tabular}{@{}p{0.45\textwidth}p{0.6\textwidth}@{}}
\toprule
\textbf{System} & \textbf{Result (as of September 2026)}\\
\midrule
Aristotle (Harmonic), Aletheia (DeepMind) & 13 of some 700 Erd\H{o}s problems resolved \cite{Erdos2026,BloomErdos}\\
AlphaEvolve (DeepMind) & progress on 23 of 67 problems \cite{GGTW2025}\\
Axiom & 4 open problems; prime-gaps lowered to$<200$\\
OpenAI & 10 open problems (announced 1 August)\\
Gemini & 6 of the 10 problems of {\it First Proof} \cite{FirstProof2026}\\
GPT-5.6, Aletheia & more than 50\% of FrontierMath Tier 4 \cite{Tier4paper}\\
\bottomrule
\end{tabular}
\smallskip
\caption{Some highlights in the first nine months of 2026.}
\label{tab:score}
\end{table}

\subsubsection{FrontierMath}
FrontierMath \cite{Glazer2024,EpochFrontierMath} is a benchmark run by Epoch AI, an organisation devoted to the measurement of AI, with funding from OpenAI. It was launched in December 2024, in direct response to the aforementioned Olympiad silver medal. Very well, the thinking went, the machine can do the cleverest of school problems; what can it do with more advanced problems? The answer, at the time, was reassuring: a mere two per cent.

Tiers 1 to 3, some three hundred problems at advanced undergraduate and beginning graduate level, were in place by March 2025. Think of the Putnam competition in place of the IMO, and then one level up. At the launch, the feeling was that this would be safe for years, and that the machines would not exceed five per cent. Within months they were at 10 to 25 \%. So the organisers scratched their heads and went to the next level.

This was Tier 4, and this is where I joined. In May 2025, thirty of us were flown to Epoch's headquarters in Berkeley, an episode which found its way into the {\it Scientific American} and the {\it Financial Times} \cite{SciAm2025}. Our task was to invent problems:  one on which a professional mathematician would work, or would be given to a collaborator, to an advanced postdoc, or to a very, very gifted graduate student.
One constraint to the benchmark made such a large-scale benchmark possible: however hard and however abstract the question, the answer must be a single exact number of at least some 20 digits. This removes any need to check a proof, which is the purpose of the FirstProof project \cite{FirstProof2026} whose batches contain 10 problems - unlike the hundreds to thousands in Tiers 1 to 4 - because human checking of proofs is very time and resource consuming.

To give the flavour, here is a public sample, due to Ken Ono (see original at \cite{EpochFrontierMath})
\KeyIdea{
Define, on the upper half-plane,
\[
F(z):=1+\sum_{n\geq1}\frac{e^{2\pi i n^2 z}}{\prod_{j=1}^{n}(1+e^{2\pi i j z})^2},\qquad
G(z):=\prod_{n\geq1}\frac{(1-e^{2\pi i n z})(1-e^{2\pi i(2n-1)z})}{1+e^{2\pi i n z}} .
\]
Let $\ell_1$ be the smallest prime such that: $-\ell_1$ is the discriminant of the ring of integers of $\Q(\sqrt{-\ell_1})$; the class number of this field is a prime $\ell_2\geq5$; the class of $\ell_2$ is a primitive root modulo $\ell_1$; and the Mordell-Weil group over $\Q$ of $Y^2=X^3-\ell_1^2X$ is $\Z/2\Z\times\Z/2\Z$. Let
\[
\alpha:=\lim_{y\to0^+}\left.\left(F-G+\frac{F}{G}-\frac{G}{F}\right)\right|_{z=\frac{\ell_1}{4\ell_2}+iy},
\]
with minimal polynomial $P_\alpha$ and splitting field $K_\alpha$. \\
Compute $\Omega:=[K_\alpha:\Q]^{-1}\left(P_\alpha(\ell_1)+P_\alpha(\ell_2)\right)^{\ell_2-1}$.
}

The answer is a single integer of some seventy digits.
Just read the problem. Unless one is an expert in analytic and algebraic number theory, the question does not so much as parse even to other mathematicians.
Ken was the most excited of all of us in Berkeley. I should imagine that seeing a problem of this kind fall is what moved him, in December 2025, to leave academia for Axiom, Carina Hong's famous start-up. He was among the first mathematicians of such prominence to leave for a private AI company; many have followed this year.

\begin{figure}[t]
\centering
\optionalimage{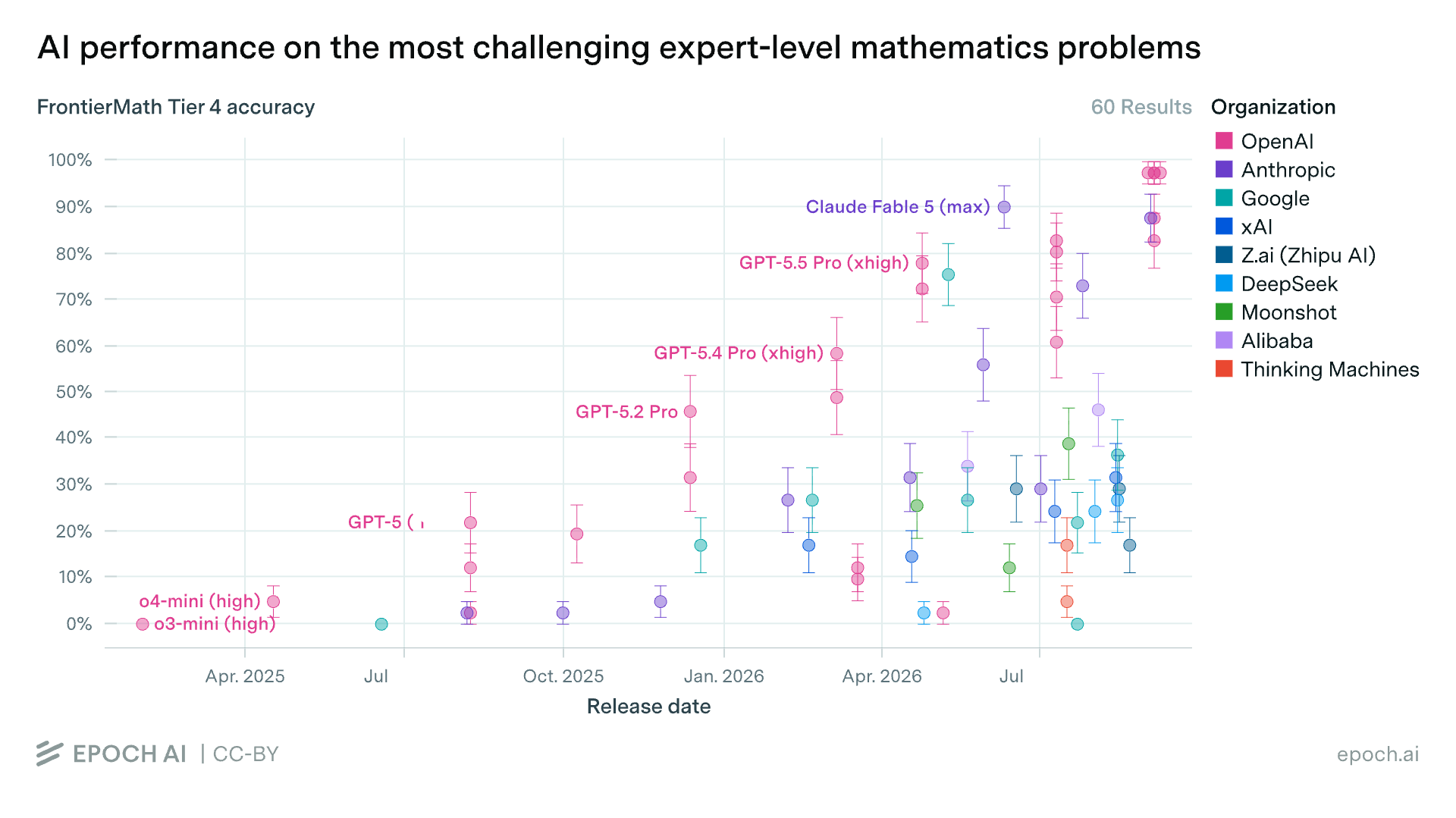}{0.9\textwidth}{2.9in}{FrontierMath Tier 4 accuracy against model release date (Epoch AI, CC-BY)}
\caption{Accuracy on FrontierMath Tier 4 by date of release of the model, from early 2025 to the autumn of 2026 \cite{EpochFrontierMath}.}
\label{fig:tier4}
\end{figure}

Figure~\ref{fig:tier4} shows what happened next. In April 2025, when we were still speaking of models with names like o3-mini, the scores were indistinguishable from zero. By October 2025 the best stood at about ten per cent. By the summer of 2026, Claude Fable had reached ninety, and the most recent models of OpenAI sit a little higher still, with those of Anthropic close behind; and these are only the public models, for each company has internal ones to which the rest of us have no access. In plain words: {\it the generic research problem in any branch of mathematics, so long as it is well defined, can be solved 90\% of the time, fully automatically by AI.}

So, in May of this year, thirty became five hundred of us: a hundred mathematicians and theoretical physicists in each of five different cities across the world, convened to write what is informally ``Tier 5''. These are {\it open} problems, to which nobody knows the answer, graded in advance by their proposers and referees as moderately interesting, solid result, major advance, or breakthrough. About fifty were accepted. The state of play at the time of the lecture, three months later, is in Table~\ref{tab:tier5}.

\begin{table}[t]
\centering
\renewcommand{\arraystretch}{1.2}
\small
\begin{tabular}{@{}lcccc@{}}
\toprule
 & \textbf{Moderately interesting} & \textbf{Solid result} & \textbf{Major advance} & \textbf{Breakthrough}\\
\midrule
Solved / posed & 6 / 22 & 2 / 18 & 1 / 6 & 1 / 3\\
\bottomrule
\end{tabular}
\smallskip
\caption{FrontierMath open problems, autumn 2026: ten of forty-nine solved, five by AI alone and five by humans with AI \cite{EpochFrontierMath}.}
\label{tab:tier5}
\end{table}

Among those solved by AI alone is a construction of Littlewood polynomials with large normalised Mahler measure. Among those solved by human and machine together are a counterexample concerning the core in approval-based committee elections, classed as a major advance, and a problem on Ap\'ery-style irrationality proofs via linear recurrences, classed as a breakthrough. 
Of particular interest to me is the construction of a polynomial whose Galois group is the Mathieu group $M_{23}$, the one sporadic group for which the inverse Galois problem had resisted all attempts. That had been open for a very long time. It no longer is.

\subsubsection{Recent striking examples}
\begin{table}[t]
\centering
\renewcommand{\arraystretch}{1.25}
\small
\begin{tabular}{@{}lp{0.2\textwidth}p{0.4\textwidth}l@{}}
\toprule
\textbf{2026} & \textbf{System} & \textbf{Result} & \textbf{Open for}\\
\midrule
May & GPT-5.4 & Erd\H{o}s unit distance conjecture: {\it counterexample} \cite{UnitDist2026} & $\sim80$ years\\
July & Claude Fable & Jacobian conjecture for $n\geq3$: {\it counterexample} & Smale's 16th\\
August & OpenAI (internal) & a non-sofic group & $\sim30$ years\\
August & Claude (internal) & $S^6$ is a complex manifold & $\sim80$ years\\
September & Claude (internal), OpenAI Astra & Navier-Stokes (to be confirmed) & Millennium\\
\bottomrule
\end{tabular}
\smallskip
\caption{Five results of the summer of 2026.}
\label{tab:striking}
\end{table}

Benchmarks are one thing. Table~\ref{tab:striking} is just a subset of 2026 so far.
In May, a counterexample was found to the unit distance conjecture which Erd\H{o}s posed in 1946. What is instructive is how: by importing methods from an adjacent part of algebraic number theory which combinatorialists had not been in the habit of using. The argument has since been improved upon by Will Sawin, who found a better proof. Humans are still in the loop, and I do not think we are about to lose our jobs. What we have, rather, is an unbelievably capable assistant at our disposal.

In July, Claude Fable found a counterexample to the Jacobian conjecture in dimension $n\geq3$. The conjecture, which is the sixteenth on Smale's list, asks whether a polynomial map of $\C^n$ with constant non-zero Jacobian determinant must have a polynomial inverse. The case $n=1$ is trivial. Now every mathematician I know is fixated on $n=2$, which is the hard one. I have myself spent a month chatting with the bots about it. Two polynomials in two variables: what on earth could be difficult? Yet so far nobody, of either species, has found a counterexample or a proof.

In August, an internal model at OpenAI produced the first example of a group which is not sofic, one of the ten results announced on the first of that month. I confess that this is above my pay grade in group theory. In the same month came the result perhaps most striking to the present audience: an internal Claude model found that the six-sphere $S^6$ admits a complex structure, settling a question which goes back to Hopf. The paper runs to some 120 pages. I believe it is now accepted as correct, and for a very particular reason: Claude then auto-formalised the entire proof into Lean. One runs the code; it compiles. This is exactly the role for Mathlib which I advertised in Section \ref{sec:topdown}, and it is, I think, how the refereeing of long proofs will be done from now on.

And then, at about five o'clock on the very morning of this lecture, came the announcement that a first Millennium Prize Problem, Navier-Stokes, has been solved. This must, of course, be auto-formalised (it has been by now) and checked by a committee (will take more than a year and rightfully so); it remains to be confirmed. I admit to some surprise at the provenance. Knowing that leading experts on the singularities of fluids, such as Javier G\'omez-Serrano, had been working with DeepMind \cite{Wang2025}, I had been saying since January that I expected DeepMind to announce a resolution of Navier-Stokes by the end of the year. Anthropic and OpenAI seem somehow to have got there first. Hence my uncertainty, in Section \ref{sec:topdown}, as to whether six Millennium Problems remain, or five.

Some comments are in order.
First, there is as yet no discernible pattern as to which {\it type} of problem is likely to fall. It is true that most of the problems fed to the machines have been in combinatorics and number theory, but that reflects the community which has been doing the feeding and the availability of curated lists, notably the seven hundred or so Erd\H{o}s problems compiled by Thomas Bloom \cite{BloomErdos}. The complex structure on $S^6$ is certainly not combinatorics.

Second, and this I find reassuring, nearly all of these results were prompted by serious mathematicians, who knew what to ask and could recognise an answer. There is one notable outlier: Erd\H{o}s problem 1196 fell to an amateur with, I am told, school-level mathematics, who somehow extracted a correct proof. But on the whole it is the professionals who are doing the asking, and we remain very much in the loop.

Third, to return to the Birch Test: the reader will have noticed that the entries of Table~\ref{tab:striking} are automatic, are as interpretable as a theorem can be, and are hardly trivial. Whether they make Birch happy, I have yet to ask him.

\section{The Mathpocalypse and Prospectus}
\begin{figure}[t]
\centering
\optionalimage{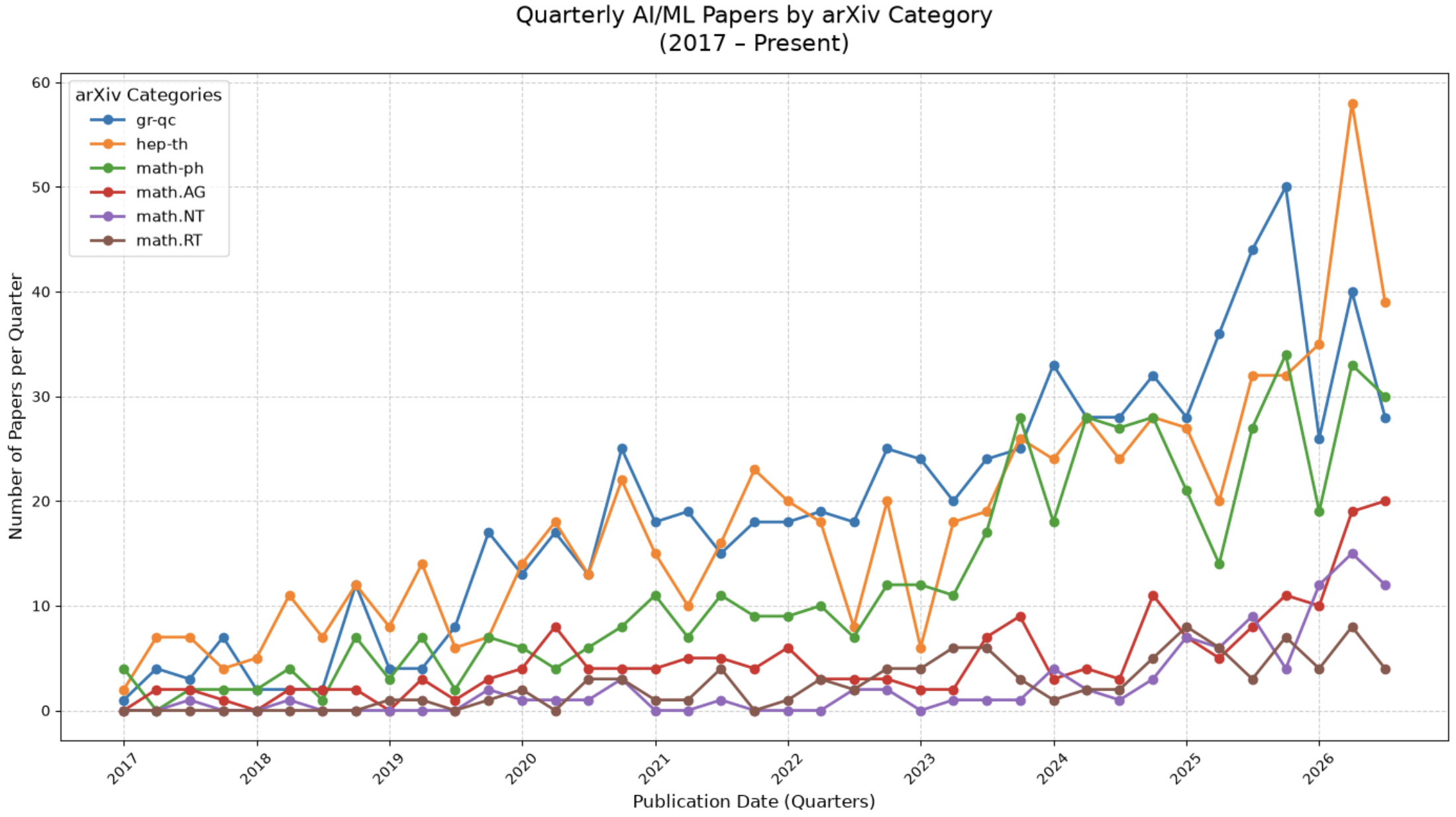}{0.85\textwidth}{2.7in}{Quarterly AI/ML papers by arXiv category, 2017 to present}
\caption{Papers per quarter with AI and machine-learning terms in title or abstract, in six arXiv categories of theoretical physics and pure mathematics.}
\label{fig:arxiv}
\end{figure}

Let me draw to a close with some sociology, for the field that is exponentiating.
Figure~\ref{fig:arxiv} is a plot which, naturally, I asked ChatGPT to make for me. It scraped the arXiv for the appearance of ``AI'', ``machine learning'', ``neural network'' and kindred terms in the titles and abstracts of papers in the unlikely categories, and counted them quarter by quarter since 2017. The physicists of \texttt{hep-th}, \texttt{gr-qc} and \texttt{math-ph} remain ahead of the pure mathematicians, with fifty or so papers a quarter, but \texttt{math.AG} is catching up quickly, at twenty a quarter and climbing, with \texttt{math.NT} and \texttt{math.RT} behind. A decade ago every one of these curves was at zero.

\begin{figure}[t]
\centering
\optionalimage{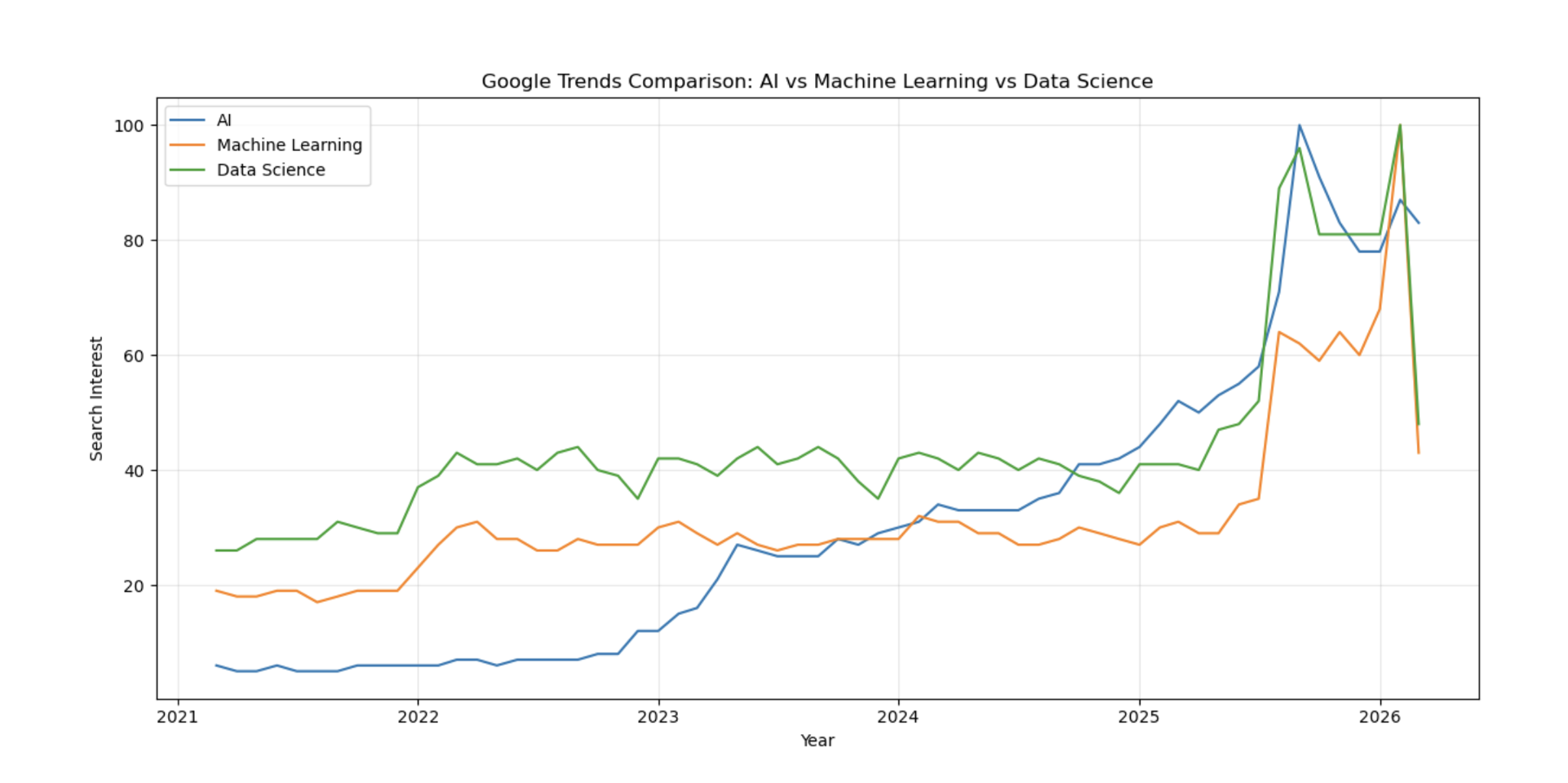}{0.8\textwidth}{2.5in}{Google Trends: AI vs Machine Learning vs Data Science, 2021 to 2026}
\caption{Search interest in ``AI'', ``machine learning'' and ``data science'' over the past five years.}
\label{fig:trends}
\end{figure}

Even the vocabulary has shifted. If one consults Google Trends (Figure~\ref{fig:trends}), the phrases ``machine learning'' and ``data science'' are steadily being absorbed into the two letters ``AI''. In 2017, undergraduates were queueing to take courses in data science, the very migration which first aroused my curiosity. That is now beside the point, and so, I am afraid, is much of machine learning as a craft, because hardly anybody codes any more. Our next generation may never quite know what coding meant, since the machine already codes better than they, or we, could.

I have had to follow the trend myself. In 2023 we launched a journal with World Scientific under the title {\it International Journal of Data Science in the Mathematical Sciences}. It was a terrible title and entirely my fault; who could possibly remember it? It has been relaunched this year as {\it AI in Mathematics and Theoretical Physics} \cite{AIMTP}, and I warmly encourage the reader to submit. I understand that Professor Yau, too, has a new journal on AI for mathematical research which is just about to appear, so that authors will be spoilt for choice.

\begin{figure}[t]
  \centering
  
  \begin{minipage}[t]{0.48\textwidth}
    \centering
    {\large \textbf{1930 math}}\\[6pt]
    \fbox{\includegraphics[width=0.92\linewidth, height=5.2cm, keepaspectratio=true]{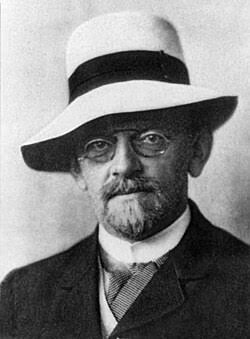}}\\[8pt]
    {\textit{``Wir m\"ussen wissen,\\ wir werden wissen.''}}\\[4pt]
    {\footnotesize David Hilbert, K\"onigsberg, 1930}
  \end{minipage}%
  \hfill
  \begin{minipage}[t]{0.48\textwidth}
    \centering
    {\large \textbf{\color{red!75!black}2026 math}}\\[6pt]
    \fbox{\includegraphics[width=0.92\linewidth, height=5.2cm, keepaspectratio=true]{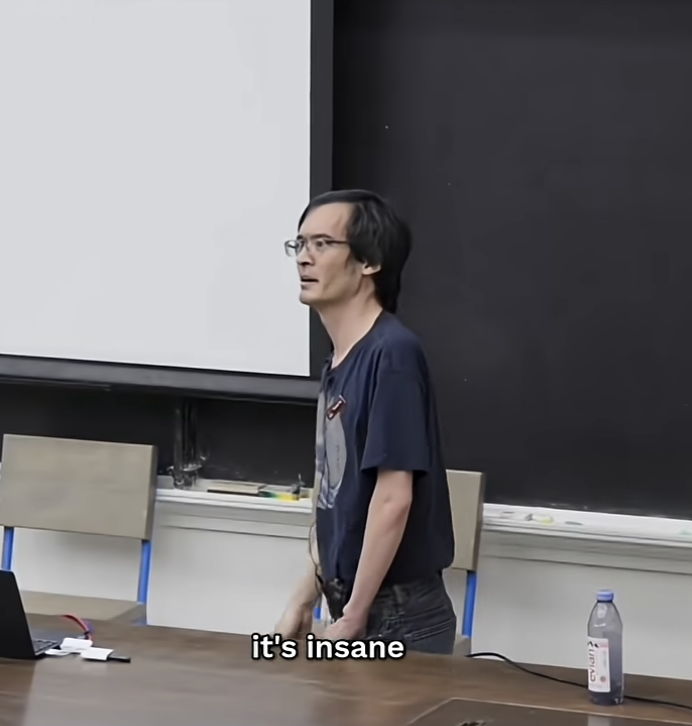}}\\[8pt]
    { \textit{``We have to\\ slow down.''}}\\[4pt]
    {\footnotesize Terence Tao, 2026}
  \end{minipage}
  \caption{Before and after Oct, 6, 2026.}
\label{f:HilbertTao}
\end{figure}

As I am writing this report for a lecture I gave a month ago, arguably the most dramatic event in the history of mathematics transpired on October 6th, 2026.
It has come to be known as the ``Mathpocalypse'' \cite{Oct6}, when OpenAI literally dumped 722 mathematical papers addressing 372 open problems across all fields of mathematics. Even the most adamant skeptic to ``AI for Maths'' was shocked \footnote{A personal anecdote. Over the years since 2017 I would regularly see Prof.~Sir Andrew Wiles over lunch at Merton College, Oxford. He has always astutely given calm counter-examples to my enthusiasm on AI abilities in mathematics. Shortly after Oct 6th, 2026, when I met him over lunch, he said, ``I think I have chosen a good time to retire.'' Meanwhile, his successor, James Maynard, is very engaged in conversations regarding the future of the field because of AI. I was honoured to share a recent podcast with him for the Economist \cite{economist}.}.
On the contrary, one of the most influential supporters has become more troubled; a slide (q.v.~Fig.~\ref{f:HilbertTao}) I have since added to my talks succinctly shows this (still taken from one of Tao's recent talks).

Over the past few months, mathematicians across the globe have been trying to digest the dump.
It is generally agreed that the results range from Fields Medal level work to claims which are simply wrong.
OpenAI - a private company with no obligation to academic rigour - did not even auto-formalize everything into Lean, let alone check through all the papers. I suppose this is what bothered the community so much. Consequently, several ``alliances \cite{AHM} for human mathematics'' have been formed to align AI with the values of the human mathematical community which took millennia to mature.

What should the working mathematician make of all this? I have tried to report rather than to predict, since any prophecy I might offer would be out of date before the ink dried, as the fate of my own reviews \cite{He2024Triumvirate,burtsev2026how} amply shows. But a few things seem to be secure. The first is that the triumvirate has become a single enterprise: the language model supplies the vibe, the Platonic data supplies the experiment, and Lean supplies the verdict. The second is that trust in mathematics is migrating from the referee to the Lean compiler, and that this is to be welcomed. The third is that a human must still decide what is worth asking. Every result in Table~\ref{tab:striking} began with mathematicians who cared about the question. The fourth is that AI has clearly demonstrated that  it can solve problems - even open ones - faster than humans, but can it create new problems, new conjectures and new fields? In a way, this is the heart of the Birch Test \cite{HeBurtsev2024}. Some of my colleagues think this is not possible since AI can never become that creative. I am more neutral on this; it already takes human creativity and ingenuity to solve open problems, there's no reason why creating new fields and visions should not be integrable from such steps already taken.

In the absolute worst case when all of discovery becomes fully automated, from raising the problem to Lean verification, to writing the paper iteratively until human digestible, there should still be value in human mathematics. Back in 2025, I joked that in the far future, should such a day come, mathematics would become more like the humanities \cite{HLF}. Indeed, scholars in English and philosophy departments have been writing and sharing commentaries on Shakespeare or Plato for centuries. That ``far future'' does not seem so far now.

It is curious to reflect upon the road we have travelled in these pages. We began with geometry: from Euler to Gauss to Riemann, and thence to Calabi and to Yau. The hyphen took us into physics, and to the great catalogues of manifolds which their work inspired. Those catalogues, forty years on from the proof, happened to be the first pure mathematics shown to a neural network, on a laptop, by a sleep-deprived father. That a single two-page announcement in the {\it Proceedings of the National Academy of Sciences} \cite{Yau1977} should have given so much, to so many disciplines, some of which did not exist when it was written, is the best argument I know for the unity of knowledge:
\[
\begin{array}{ll}
\includegraphics[trim=0mm 0mm 0mm 0mm, clip, angle=0, width=2in]{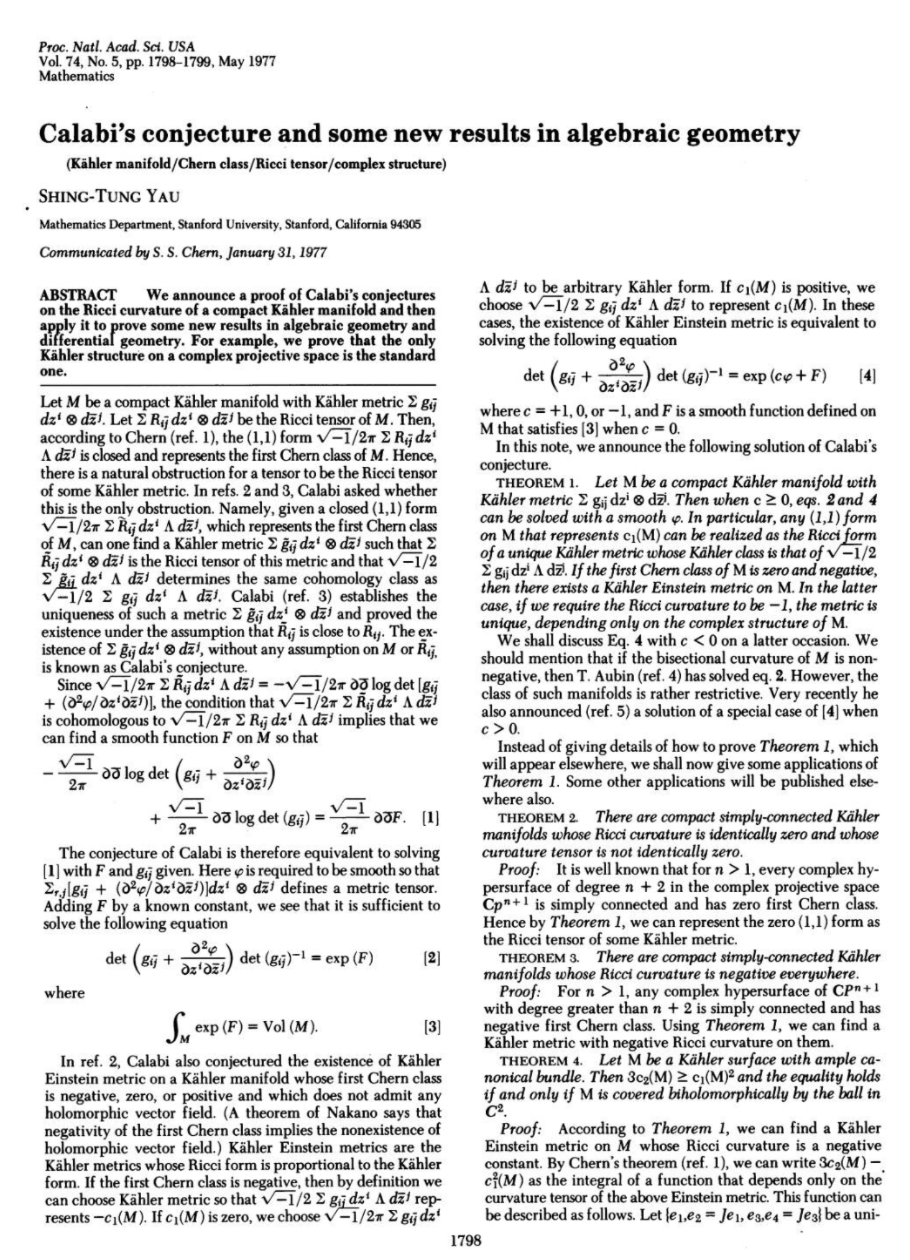} &
\includegraphics[trim=-20mm 0mm 0mm 0mm, clip, angle=0, width=4in]{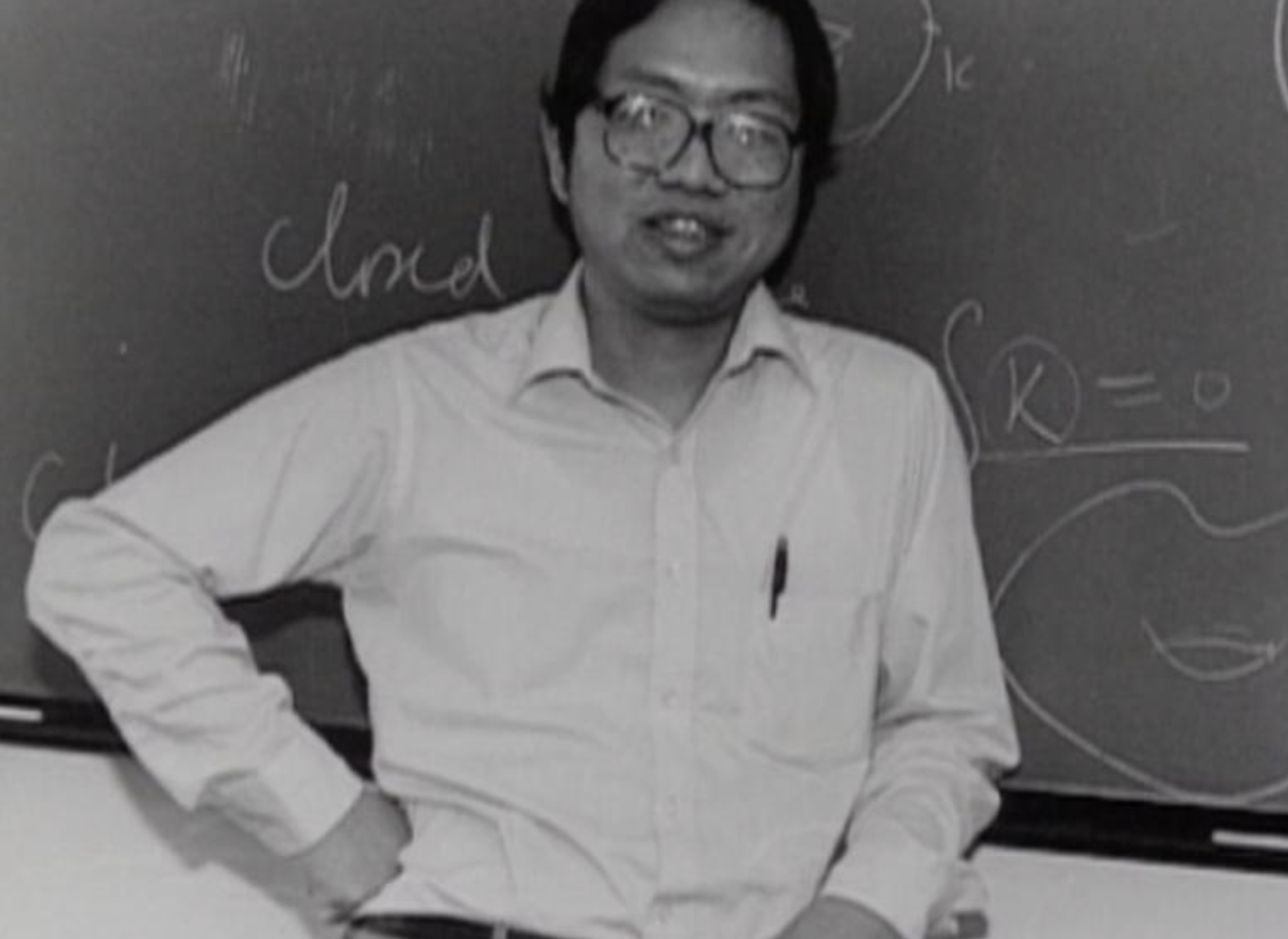}
\end{array}
\]

So, to that paper and to its author, and to all the friends, old and new, gathered in the autumnal rain of this jubilee: a very happy fiftieth.


\end{document}